\documentclass[12pt]{amsart} 
\usepackage{amssymb}

\usepackage{euscript}

\let\cal\mathcal

\makeatletter \@mparswitchfalse \makeatother

\newtheorem{theorem}{Theorem}
\newtheorem{lemma}{Lemma}
\newtheorem{corollary}{Corollary}
\newtheorem{proposition}{Proposition}
\newtheorem{remark}{Remark}

\newtheorem{definition}{Definition}
\newtheorem{claim}{Claim}

\newcommand{\pprec}{\prec\!\prec}

\def\eqref#1{(\ref{eq#1})}
\def\eqlabel#1{\label{eq#1}}

\def\N{\mathbb N}

\def\R{\mathbb R}
\def\M{\cal{M}}
\let\<\langle

\def\T{\tau}

\numberwithin{equation}{section}

\begin{document}


\let\\\cr
\let\union\bigcup
\let\inter\bigcap
\def\supp{\operatorname{supp}}
\def\sup{\operatorname{sup}}
\def\inf{\operatorname{inf}}
\def\Im{\operatorname{Im}}
\def\dim{\operatorname{dim}}
\def\Span{\operatorname{span}}
\def\cord{\operatorname{cord}}
\def\lim{\operatorname{lim}}
\def\Re{\operatorname{Re}}
\def\sqn{\operatorname{sqn}}
\def\log{\operatorname{log}}
\def\max{\operatorname{max}}
\def\min{\operatorname{min}}
\let\emptyset\varnothing
\def\exp{\operatorname{exp}}

\title{Conjugate Operators for Finite Maximal Subdiagonal Algebras}
\author{Narcisse Randrianantoanina}
\address{Depatment of Mathematics, University of Texas, Austin, TX 78712}
\curraddr{Mathematical Sciences Research Institute\\
1000 Centennial Drive\\
Berkeley, CA 94720}
\email{nrandri@math.utexas.edu}
\subjclass{46L50, 46E15; Secondary: 43A15,  47D15}
\keywords{von-Neumann algebras, conjugate functions, Hardy spaces}
\thanks{The  research reported here was undertaken during the author's visit
at the Mathematical Sciences Research Institute (MSRI) in Berkeley, CA. 
Research at the MSRI is supported in part by the NSF grant DMS-9022140  }

\begin{abstract}
Let $\M$ be a von Neumann algebra with a faithful normal trace $\T$,
and let $H^\infty$ be a finite, maximal, subdiagonal algebra of $\M$.
Fundamental theorems on conjugate functions for weak$^*$\!-Dirichlet
algebras are shown to be valid for non-commutative $H^\infty$. In
particular the conjugation operator is shown to be a bounded linear map 
from $L^p(\M, \T)$ into $L^p(\M, \T)$ for $1 < p < \infty$, and to be a
continuous map from  $L^1(\M,\T)$ into $L^{1, \infty}(\M,\T)$. We also 
obtain that if an operator $a$ is such that $|a|\log^+|a| \in L^1(\M,\T)$
then its conjugate belongs to $L^1(\M,\T)$.  Finally, we present some partial
extensions of the classical Szeg\"o's theorem to the non-commutative setting.
\end{abstract}

\maketitle

\section{introduction}
The theory of conjugate functions has been a strong motivating force behind
various aspects of harmonic analysis and abstract analytic function spaces.
This theory which was originally developed for functions in the circle group
$\mathbb{T}$ has found many generalizations to  more abstract  settings
 such as 
Dirichlet algebras in \cite{DEV} and weak*-Dirichlet algebras in \cite{HR}.
Results from this theory
 have been proven to be very fruitful for studying Banach space
properties of the Hardy spaces (and their relatives) associated with the
 algebra involved (see for instance \cite{BO6} and \cite{LAN}).
 
Let $\M$ be a von Neumann algebra with a faithful, normal finite trace $\T$. 
 Arveson introduced in \cite{AV}, as non-commutative
analogues of weak$^*$\!-Dirichlet algebras, the notion of finite, maximal 
subdiagonal algebras of $\M$ (see definition below).
 Subsequently several authors
studied the (non-commutative) $H^p$-spaces associated with  such algebras
(\cite{kT}, \cite{MAR}, \cite{MMS}, \cite{SAI2}, \cite{SAI1}).
 In \cite{MAR},  the notion of
harmonic conjugates  was introduced 
for maximal subdiagonal algebras
 generalizing the notion of conjugate functions for 
weak*-Dirichlet algebras and it was proved that the operation
of conjugation is bounded in $L^2(\M, \T)$.

The main objective of this paper  is to
combine the spirit of \cite{HR} with that of
\cite{MAR} to get a more constructive
 definition of conjugate operators for the setting of non-commutative
maximal subdiagonal algebras;
and  to study different properties of conjugations
 for these non-commutative settings.  We prove that most fundamental
theorems on conjugate operation on Hardy spaces associated with
weak$^*$\!-Dirichlet (see \cite{DEV} and \cite{HR})
 remain valid for Hardy spaces associated with finite subdiagonal
algebras. In particular, we show that the conjugation operator is a
bounded map from $L^p(\M,\T)$ into $L^p(\M, \T)$ for $1< p < \infty$,
 and from $L^1(\M, \T)$ into $L^{1, \infty}(\M, \T)$.
 We conlude that, as in commutative case, 
 (non-commutative) $H^p$ is a complemented subspace of  $L^p(\M, \T)$
for $1<p<\infty$. 

Many results in harmonic analysis can be   
deduced from the classical Szeg\"o's 
theorem. This very classical fact, although is valid  for the more abstract
setting of weak*-Dirichlet algebras, is still unknown for the non-commutative
 case. The last part of this paper is devoted to various results related to
Szeg\"o's theorem.

We refer to \cite{N}, \cite{SEG} and  \cite{TAK} for general information
 concerning
von Neumann algebras as well as basic notions of non-commutative
integration,  to \cite{D1} and \cite{LT} for Banach space theory 
and to \cite{HEL} and \cite{ZYG} for basic definitions from
harmonic analysis. %
\section{Definitions and preliminary results}

Throughout, $H$ will denote a Hilbert space and $\M \subseteq
\cal{L}(H)$ a von Neumann algebra with a normal, faithful finite trace
$\T$. A closed densely defined operator $a$ in $H$ is said to be {\em
affiliated with} $\M$ if $u^* au = a$ for all unitary $u$ in the
commutant $\M'$ of $\M$. If $a$ is a densely defined self-adjoint
operator on $H$, and if $a = \int^\infty_{- \infty} s d e^a_s$ is its
spectral decomposition, then for any Borel subset $B \subseteq \R$,
we denote by $\chi_B(a)$ the corresponding spectral projection
$\int^\infty_{- \infty} \chi_B(s) d e^a_s$. A closed densely defined
operator on $H$ affiliated with $\M$ is said to be {\em
$\T$-measurable} if there exists a number $s \geq 0$ such 
that $\T(\chi_{(s, \infty)} (|a|)) < \infty$.

The set of all $\T$-measurable operators will be denoted by
$\overline{\M}$. The set $\overline{\M}$ is a $*$\!-algebra with respect to the
strong sum, the strong product, and the adjoint operation \cite{N}. %
For $x \in \overline{\M}$, the generalized singular value function $\mu (x)$
of $x$ is defined by
$$
\mu_t(x) = \inf \{ s \geq 0: \T(\chi_{(s, \infty)} (|x|)) \leq t \},
\quad \text{ for } t \geq 0.
$$
The function $t \to \mu_t(x)$ from $(0, \T(I))$ to $[0, \infty)$
is right continuous, non-increasing and is the inverse of the
distribution function $\lambda (x)$, where $\lambda_s(x) = \T(\chi_{(s,
\infty)}(|x|))$, for $s \geq 0$. For a complete study of $\mu(.)$ and
$\lambda(.)$ we refer to \cite{FK}.

\begin{definition}
Let $E$ be an order continuous rearrangement invariant (quasi-)
Banach function space on $(0, \T(I))$. We define the symmetric space
$E(\M, \T)$ of measurable operators by setting:
\begin{align*}
E(\M, \T) &= \{ x \in
\overline{\M}\quad ; \quad \mu(x) \in E \} \quad \text{and} \\
\|x\|_{E(\M,\T)} &= \| \mu(x)\|_E, \text{ for } x \in E(\M,\T).
\end{align*}
\end{definition}

It is well known that $E(\M, \T)$ is a Banach space
(resp. quasi-Banach space) if $E$ is a Banach
space (resp. quasi-Banach space),
 and that if $E = L^p(0, \T(I))$, for $0 < p < \infty$, then 
$E(\M, \T)$ coincides with the usual non-commutative $L^p$\!-space
associated with $(\M, \T)$. We refer to \cite{CS}, \cite{DDP1} and \cite{X}
 for more detailed discussions about these
spaces. For simplicity we will always assume that the trace $\tau$ is 
normalized. %

The following definition isolates the main topic of this paper.

\begin{definition}
Let $H^\infty$ be a weak$^*$\!-closed unital subalgebra of $\M$ and let
$\Phi$ be a faithful, normal expectation from $\M$ onto the diagonal $D
= H^\infty \cap (H^\infty)^*$, where $(H^\infty)^* = \{ x^*,\  x \in
H^\infty \} $. Then $H^\infty$ is called a finite, maximal,
subdiagonal algebra in $\M$ with respect to $\Phi $ and $\T$ if:
\begin{itemize} 
\item[(1)] $H^\infty + (H^\infty)^*$ is weak$^*$\!-dense in $\M$;
\item[(2)] $\Phi (ab) = \Phi(a) \Phi(b)$ for all $a,b \in H^\infty$;
\item[(3)] $H^\infty$ is maximal among those subalgebras satisfying (1)
and (2);
\item[(4)] $\T \circ \Phi = \T$.
\end{itemize}
\end{definition}

For $0 < p < \infty$, the closure of $H^\infty$ in $L^p(\M, \T)$ is
denoted by $H^p(\M,\T)$ (or simply $H^p$) and is called the Hardy space
 associated with the subdiagonal algebra $H^\infty$. Similarly, 
  the closure of $H^\infty_0 = \{ x \in H^\infty;
\Phi(x) = 0 \}$ is denoted by $H^p_0$.

Note that $\Phi$ extends to $L^2(\M, \T)$ and this extension is an
orthogonal projection from $L^2(\M, \T)$ onto $[D]_2$, the closure of
$D$ in $L^2(\M, \T)$. Similarly, since $ \| \Phi(x)\|_1 \leq \|x\|_1$
for every $x \in \M$, the operator $\Phi$ extends 
uniquely to a projection of norm one from  $L^1(\M, \T)$ onto
$[D]_1$, the closure of $D$ in $L^1(\M, \T)$.

\section{Harmonic conjugates}

Let $\cal{A} = H^\infty + (H^\infty)^*$. Since $\cal{A}$ is
weak$^*$\!-dense in $\M$, it
is norm dense in $L^p(\M, \T)$, where $ 1 \leq p < \infty$.

Note that $H^\infty$ and $(H^\infty_0)^*$ are orthogonal in $L^2(\M,\T)$.
This fact implies that $L^2(\M,\T)=H^2 \oplus (H^2_0)^*$, and hence
that $L^2(\M,\T)= H^2_0 \oplus (H^2_0)^* \oplus [D]_2$.

Let $a \in \cal{A}$. Then $a$ can be written as  $a_1 + a_2^* + d$ where
$a_1$ and   $a_2$ belong to $ H^\infty_0$ and $d \in D$.
 In fact, $ a = b_1 + b_2^*$
with $b_1, b_2 \in H^\infty $ and set $d = \Phi(b_1) + \Phi(b_2^*) \in
D$ and $a_i = b_i - \Phi(b_i)$, for $i = 1,2$.

Since $H^\infty_0$ and $(H^\infty_0)^*$ are orthogonal subsets of
$L^2(\M, \T)$, this decomposition is unique. For $u = u_1 + u_2^* + d$ in
$\cal{A}$ , we define $\tilde{u} = i u_2^* - i u_1$. Then $\tilde{u} \in \M$
and $u + i \tilde{u} = 2 u_1 + d \in H^\infty$. 
The operator $\tilde{u}$ will be called the {\bf conjugate} of $u$.

 Our main
goal is to study the operation that takes $u \in \cal{A}$ into $\tilde{u} \in
\M$ as linear operator between non-commutative $L^p$-spaces.
 In particular we will 
 extend ``$\sim$'' to $L^p(\M, \T)$ for $1 \leq p <
\infty$. It should be noted that if $\M$ is commutative, then the above
definition coincides with the the definition of conjugate functions for
weak$^*$\!-Dirichlet algebras studied in \cite{HR}.

\begin{remark} 
\begin{itemize}
\item[(i)] $\Phi(\tilde{u}) = 0$, for all $u \in \cal{A}$.
\item[(ii)] If $u = u^*$, then the uniqueness of the decomposition
implies that $u_1 = u_2$ and $d=d^*$.
 Therefore  if $u=u^*$ then $\tilde{u} = \tilde{u}^*$.
\item[(iii)] For $u = u_1 + u_2^* + d \in \cal{A}$ and 
$\tilde{u} = i(u_2^* - u_1)$, the above observation implies that
$u_2^* \perp u_1 $ in $L^2(\M,\T)$, so
$$
\| \tilde{u} \|^2_2 = \| u_2^* - u_1 \|^2_2 = \| u_2^* \|^2_2 + \| u_1
\|^2_2,
$$
and since $L^2(\M,\T) = H^2_0 \oplus (H^2_0)^* \oplus [D]_2$ we get,
 $$\| u\|^2_2 = \| u_1 \|^2_2 + \| u_2^* \|^2_2 + \| d \|^2_2,$$
 which implies
that $\| \tilde{u} \|_2 \leq \| u \|_2$.
\end{itemize}
\end{remark}
As a consequence of (iii), we get the following theorem:
\begin{theorem}
There is a unique continuous linear map ``$\sim$'' from
$L^2(\M, \T)$ into $L^2(\M,\T)$ that coincides with ``$\sim$'' in $\cal{A}$.
This map is of norm 1, and if $u \in L^2(\M, \T)$ then $ u + i \tilde{u}
\in H^2$.
\end{theorem}
We remark that Marsalli has recently proved a version of Theorem~1 (see
\cite{MAR} Corollary~10):
 he showed that the conjugation operator is bounded in $L^2(\M,\T)$
with bound less than or equal to $\sqrt{2}$.
\subsection*{Extension of the operator ``$\sim$''
 to  $L^p(\M,\T)$, $1<p<\infty$}

\

In this section, we will  extend Theorem 1 from $p=2$ to all $p$
with $1 < p < \infty$.

The following elementary lemma will be used in the sequel; we
will include its  proof for completeness.

\begin{lemma}
Let $m \in \N$ and $a_1,a_2, \ldots, a_m \in \overline{\M}$. If
$\frac{1}{p_1} + \frac{1}{p_2} + \ldots + \frac{1}{p_m}=1$, and $a_j \in
L^{p_j}(\M,\T)$ for each $ j \leq m$, then
$$
|\T(a_1 a_2 \ldots a_m)| \leq \Pi^m_{j=1} \| a_j \|_{p_j}.
$$
\end{lemma}

\begin{proof}
 Recall that, for $a, b \in \overline{\M}$, the operator  $a$ is said to be 
submajorided by $b$ and write $a \pprec b$ if
$$
\int^\alpha_0 \mu_t (a) dt \leq \int^\alpha_0 \mu_t (b) dt, \quad
\text{ for all } \alpha \geq 0.
$$
The lemma will  be proved inductively on $m \in \N$:

 For $m=2$, it is the usual
H\"older's  inequality.

Let $\frac{1}{p} + \frac{1}{q} = \frac{1}{r}$ and $b, c \in
\overline{\M}$. Then $\| bc \|_r \leq \| b \|_p \cdot \| c \|_q$: this is
a consequence of the fact \cite[Theorem 4.2\,(iii)]{FK} that
$
\mu_{(.)}^l (bc) \pprec \mu_{(.)}^l(b) \mu_{(.)}^l (c) \quad \text{
for all }l \in \N$.
So
$$
\| bc \|^r_r = \int \mu^r_t (bc)dt \leq \int \mu^r_t (b) \mu^r_t (c)
dt;
$$
then apply the usual H\"older's inequality for functions.

 Now assume
that the lemma is valid for $m=1, 2, \ldots, k$. Let $a_1, a_2,
\ldots, a_k, a_{k+1} \in \overline{\M}$ and $\frac{1}{p_1} +
\frac{1}{p_2} + \cdots + \frac{1}{p_k} + \frac{1}{p_{k+1}}
= 1$. Choose $q$ such that
$\frac{1}{q}=\frac{1}{p_k}+\frac{1}{p_{k+1}}$.  Then
\begin{align*}
\T(a_1\ldots a_{k-1} \cdot(a_k a_{k+1})) & \leq
\Pi^{k-1}_{i=1} \| a_j \|_{p_j} \cdot \| a_k
a_{k+1}\|_q \\
& \leq \Pi^{k+1}_{i=1} \| a_j \|_{p_j}.
\end{align*}
The proof is complete.
\end{proof}
\begin{theorem}
For each $1 < p < \infty$, there is a unique continuous linear
extension of ``$\sim$'' from $L^p(\M,\T)$ into $ L^p(\M,\T)$ with the
property that $f + i \tilde{f} \in H^p$ for all $f \in L^p(\M,\T)$.
Moreover there is a constant $C_p$ such that 
$$
\| \tilde{f}\|_p \leq C_p  \| f \|_p
\quad \text{ for all } f \in L^p(\M,\T).
$$
\end{theorem}

\begin{proof}
Our proof follows Devinatz's argument (\cite{DEV}) for Dirichlet algebras, 
but at number of points, certain
non-trivial adjustments have to be made to fit the non-commutative setting.

Let $u \in \cal{A}$ be  nonzero and self-adjoint; $\tilde{u}$ is self-adjoint.
Let $g=u+i \tilde{u} \in H^\infty$. Since $u = u^*$, it is of the form
$u = a+a^* + d$, where $a \in  H^\infty_0$ and $d=d^* \in D$. Recall that
$\tilde{ u}= i(a^* - a)$ so $g = 2a + d \in H^\infty$. We get that
$$
\Phi(g^{2k}) = \Phi ((2a+d)^{2 k}) = [2 \Phi(a) +
\Phi(d)]^{2 k } = \Phi (d)^{2 k}.
$$
So $\Phi((u+i \tilde{u})^{2 k}) = \Phi(d)^{2 k}$ and %
taking the adjoint, $\Phi((u - i \tilde{u})^{2 k}) = \Phi (d)^{2
k}$. Adding these two equalities, we get 
\begin{equation}
\eqlabel{3}
\Phi [ ( u + i \tilde{u})^{2 k} + (u - i \tilde{u})^{2 k}] =
2 \Phi(d)^{2 k}.
\end{equation}
Now we will expand the operators $(u + i \tilde{u})^{2 k}$ and $(u - i
\tilde{u})^{2 k}$. Note that $u$ and $\tilde{u}$  do not
necessarily commute. 

For $2\leq m\leq 2k$, let
 $\cal{S}_m = \{(r_1, r_2, \ldots , r_{m}) \in \{ 1, \ldots, 2k-1\}^{m };
\quad \sum^{m }_{j=1} r_j = 2k\}$ 
 and set $\cal{S}=\cup_{2\leq m \leq 2k} \cal{S}_m$.
For a finite sequence of integers $r=(r_1, r_2, \ldots, r_m)$, 
we set $s(r)=\sum_{j=1}^{[m/2]} r_{2j}$. Then
\begin{align*}
(u + i \tilde{u})^{2k} &= u^{2k} +(i\tilde{u})^{2k} +
 \sum_{(r_1,  \ldots, r_{m})\in
\cal{S}} (u^{r_1} (i\tilde{u})^{r_2}\ldots) +
  ((i\tilde{u})^{r_1}u^{r_2} \ldots) \\
&= u^{2k} +(i)^{2k}\tilde{u}^{2k} +
\sum_{(r_1, \ldots, r_{m}) \in \cal{S}}(i)^{s(r)}
(u^{r_1} \tilde{u}^{r_2}\ldots)  +
 (i)^{2k - s(r)} (\tilde{u}^{r_1}u^{r_2}\ldots).
\end{align*}
Similarly,
$$
(u - i \tilde{u} )^{2k} = u^{2k} +(-i)^{2k}\tilde{u}^{2k} +
\sum_{(r_1, \ldots , r_{m})
\in \cal{S}}(-i)^{s(r)} (u^{r_1} \tilde{u}^{r_2}\ldots)+
 (-i)^{2k - s(r)} (\tilde{u}^{r_1}u^{r_2}\ldots).
$$
If $\cal{K} = \{r= (r_1,r_2,\ldots,r_{2m})\in
 \cal{S}; \ s(r) \in 2\N \}$, then
$$
(u + i \tilde{u})^{2k} + (u - i \tilde{u})^{2k} = 2u^{2k} +
 2(i)^{2k}\tilde{u}^{2k} + 2 \sum_{r\in \cal{K}}
 (i)^{s(r)} (u^{r_1}\tilde{u}^{r_2}\ldots) +
  (i)^{2k - s(r)} (\tilde{u}^{r_1}u^{r_2}\ldots),
$$
so from \eqref{3}, we get
$$
\Phi(d)^{2k}= 
 \Phi(u^{2k}) + (i)^{2k}\Phi(\tilde{u}^{2k}) +
 \sum_{r \in \cal{K}}
(i)^{s(r)}  \Phi (u^{r_1} \tilde{u}^{r_2}\ldots) +
  (i)^{2k -s(r)} \Phi(\tilde{u}^{r_1}u^{r_2}\ldots).
$$

This implies 
$$
(i)^{2k} \Phi(\tilde{u}^{2k}) = \Phi(d)^{2k} - \Phi(u^{2k}) -
 \sum_{r \in \cal{K}}(i)^{s(r)}
\Phi(u^{r_1}\tilde{u}^{r_2}\ldots) +
 i^{2k -s(r)} \Phi(\tilde{u}^{r_1}u^{r_2}\ldots).
$$
Taking the trace on both sides,
$$
|\T(\tilde{u}^{2k})| \leq |\T(d^{2k})| + |\T(u^{2k})| +
 \sum_{r  \in \cal{K}} | \T(u^{r_1} \tilde{u}^{r_2}\ldots)|+
 |\T(\tilde{u}^{r_1}u^{r_2}\ldots)|.
$$
Applying  Lemma~1,  with $\frac{1}{p_j} = \frac{r_j}{2k}$, for every
$r=(r_1, r_2, \ldots, r_{m}) \in \cal{K}$, we get
$$
|\T(\tilde{u}^{2k})| \leq |\T(d^{2k})| + \| u \|^{2k}_{2k} +
 \sum_{r  \in \cal{K}} (\| u \|^{r_1}_{2k} \| \tilde{u} \|^{r_2}_{2k}
\ldots) + (\| \tilde{u} \|^{r_1}_{2k} \| u \|^{r_2}_{2k}\ldots).
$$
We observe that by the definition of $\cal{K}$,
$$ \| u \|^{2k}_{2k} +
\sum_{r \in \cal{K}} (\| u \|^{r_1}_{2k} \|
\tilde{u} \|^{r_2}_{2k} \ldots) +
 ( \| \tilde{u}\|^{r_1}_{2k} \| u \|^{r_2}_{2k}\ldots)
$$
is  equal to  the sum of the terms of the expansion of $( \| u
\|_{2k} + \| \tilde{u} \|_{2k})^{2k}$ with $\| \tilde{u} \|_{2k}$ of
even exponents between $2$ and $2k-2$, i.e., 
\begin{multline*} 
 \| u \|^{2k}_{2k} + \sum_{r
 \in \cal{K}} (\| u \|^{r_1}_{2k} \| \tilde{u} \|^{r_2}_{2k}
\ldots) + ( \| \tilde{u}\|^{r_1}_{2k}\| u \|^{r_2}_{2k} \ldots)\\
=\binom{2k}{0} \| u \|^{2k}_{2k} + \binom{2k}{2} \| u \|^{2k-2}_{2k}
\| \tilde{u} \|^{2}_{2k} 
 + \cdots + \binom{2k}{2k-2} \| u \|^{2}_{2k} \| \tilde{u}
\|^{2k-2}_{2k}.
\end{multline*}
Since $\tilde{u}$ is self-adjoint, $\T(\tilde{u}^{2k})=
\| \tilde{u}\|_{2k}^{2k}$
and hence,
$$
\| \tilde{u} \|^{2k}_{2k} \leq \| d \|^{2k}_{2k} + \| u \|^{2k}_{2k} +
\binom{2k}{2} \| u \|^{2}_{2k} \| \tilde{u} \|^{2(k-1)}_{2k} + \cdots
+ \binom{2k}{2k-2} \| u \|^{2(k-1)}_{2k} \| \tilde{u} \|^{2}_{2k};
$$
and since $\| d \|^{2k}_{2k} \leq \| u \|^{2k}_{2k}$, we have
$$
\| \tilde{u} \|^{2k}_{2k} \leq 2 \| u \|^{2k}_{2k} + \binom{2k}{2} \|
u \|^2_{2k} \| \tilde{u} \|^{2(k-1)}_{2k} + \cdots + \binom{2k}{2k-2}
\| u \|^{2(k-1)}_{2k} \| u \|^2_{2k}.
$$

Divide both sides by $\|u \|^{2k}_{2k}$ and set $X_0 = \| \tilde{u} \|_{2k}/\|
u \|_{2k}$, we have
$$
X^{2k}_{0} - \binom{2k}{2} X^{2(k-1)}_{0} - \binom{2k}{4}
X^{2(k-2)}_{0} - \cdots - 2 \leq 0.
$$
Hence, $X_0$ is less than or equal to the largest real root of the
polynomial equation
$$
X^{2k} - \binom{2k}{2}X^{2(k-1)} - \binom{2k}{4}X^{2(k-2)} - \cdots -
2= 0.
$$
If the largest root is $K_{2k}$, we have
$$
\| \tilde{u} \|_{2k} \leq K_{2k} \| u \|_{2k}.
$$
 Using Minkowski's inequality, we conclude that for every $ f \in \cal{A}$ 
(not necessarily self adjoint), we have
$$
\| \tilde{f} \|_{2k} \leq 2K_{2k} \| f \|_{2k}.
$$
Since $\cal{A}$ is dense in $L^{2k}(\M,\T)$, 
the inequality above shows that ``$\sim$'' can be
extended as a bounded linear operator from $L^{2k}(\M,\T)$ into
 $L^{2k}(\M,\T)$, so the theorem is proved for $p$ even.

For the general case, let 
$ 2 \leq p < \infty$. Choose  $k$ such that $2k \leq p
\leq 2k + 2$. By \cite{DDP2} (Theorem~2.3), $L^p(\M,\T)$ can be 
realized as a complex 
interpolation of the pair $(L^{2k}(\M,\T), L^{2k+2}(\M,\T))$, and we conclude
that ``$\sim$'' is also bounded from $L^p(\M,\T)$ into $L^p(\M,\T)$.

For $1 < p < 2$,  from the above case,``$\sim$'' is bounded from $L^{q}(\M,\T)$
 into
$L^{q}(\M,\T)$, where $\frac{1}{p} + \frac{1}{q} = 1$, and we claim that as in
 the commutative case,  $(\sim)^* = - (\sim)$.

To see this,
let  $u$ and  $v$ be self adjoint elements of   $ \cal{A}$; we have
$$
\Phi((u + i \tilde{u})(v + i \tilde{v})) = \Phi(u+i \tilde{u}) \cdot
\Phi(v + i \tilde{v}) = \Phi(u)\Phi(v),
$$
which implies that
$$
\Phi(uv + iu \tilde{v} + i \tilde{u} v - \tilde{u} \tilde{v}) =
\Phi(uv - \tilde{u} \tilde{v}) + i \Phi(u \tilde{v} + \tilde{u} v) =
\Phi(u)\Phi(v),
$$
so
$$
\T (uv - \tilde{u} \tilde{v}) + i \T(u \tilde{v} + \tilde{u} v) =
\T(\Phi(u)\Phi(v)).
$$
Since $\Phi(u)$ and $\Phi(v)$ are self-adjoint, $\T(\Phi(u) \Phi(v))
\in \R$, and also $\T (uv - \tilde{u} \tilde{v})$ and $\T(u \tilde{v}
+ \tilde{u} v) \in \R$. This implies $ \T(u \tilde{v} + \tilde{u}v)=0$
 and $ \T(u \tilde{v}) = - \T(\tilde{u} v)$. The proof is complete.
\end{proof}

\subsection*{Extension of the operator ``$\sim$'' to $L^1(\M,\T)$.}

\

The following lemma is probably known but we will include its proof for the
convenience of the reader.

\begin{lemma}
For $u \in \M$, $u \geq 0$, let $f = u + i \tilde{u}$ and $ 0 < \varepsilon
< 1$. 
\begin{itemize}
\item[(1)] $I + \varepsilon f$ has bounded inverse with $\|(I + \varepsilon
f)^{-1} \| \leq 1$.
\item[(2)] $f_{\varepsilon} = (\varepsilon I + f) (I + \varepsilon
f)^{-1} \in H^\infty$.
\item[(3)]$\Re\,( f_{\varepsilon}) \geq \varepsilon I$.
\item[(4)] $\lim_{\varepsilon \to 0} \| f_{\varepsilon} -  f \|_p =0$
 \ ($ 1\leq p <\infty$).
\end{itemize}
\end{lemma}

\begin{proof}
(1) Note that $f$ is densely defined and that, for every $x \in
D(f)$,
$$
\langle(I + \varepsilon f) x, \, x \rangle = \langle(I+ \varepsilon
u)x, \, x \rangle + i \langle \tilde{u} x, \, x \rangle.
$$
Thus $|\langle(I+ \varepsilon f) x, \, x \rangle| \geq \|x\|^2$,
which implies 
$$
\|(I + \varepsilon f) x \| \geq \|x\| \quad \text{ for all } x \in
D(f).
$$
So $I + \varepsilon f$ is one-to-one; now for every $y \in R(I +
\varepsilon f)$ define $(I + \varepsilon f)^{-1} y$ to be the
unique element such that $(I + \varepsilon f)((I+ \varepsilon
f)^{-1} y) = y$. Then $(I + \varepsilon f)^{-1}$ is linear
and, for every $y \in R(I+ \varepsilon f)$, we have $\| (I +
\varepsilon f)^{-1} y \| \leq \|y\| $. We claim  that
$R(I + \varepsilon f)$ is dense in $H$. For this, note that (using similar
estimate), $I + \varepsilon f^*$ is one-to-one;
 if $z \perp R(I +\varepsilon f)$ then
 $x \to \langle z, (I +\varepsilon f)x\rangle =0$ is continuous so
 $z \in D(I + \varepsilon f^*)$ with
 $\langle (I +\varepsilon f^*)z, x\rangle=0$ for every
 $x \in D(I +\varepsilon f)$ so $z=0$. 
Hence $(I +\varepsilon f)^{-1}$ can be  extended as  a bounded operator of
 norm $\leq 1$.
\smallskip

(2) To prove that $f_{\varepsilon} \in H^\infty$, note that
$(I + \varepsilon f)^{-1} \in \M$ with inverse $(I + \varepsilon f)
\in H^2$. %
In particular, the inverse of $(I+ \varepsilon f)^{-1}$ lies in
$L^2(\M,\T)$, so from Proposition~1.2 of \cite{MMS} (see also Proposition~1 of
\cite{SAI1}),
 there exists a unitary operator $a \in \M$ and an operator
$b \in H^\infty$ such that $(I + \varepsilon f)^{-1} = ab$. Thus
$a^*(I+\varepsilon f)^{-1} = b \in H^\infty$, and since $1 +
\varepsilon f \in H^2$, we have $a^* = b(I+ \varepsilon f) \in H^2$,
so $a^* \in H^2 \cap \M$, which implies $a^* \in H^\infty$. Now $a$ is
unitary, that is, $aa^* = a^*a = I$ (i.e $ a^* = a^{-1}$), and
$H^\infty$ is a Banach algebra, so $a = (a^*)^{-1} \in H^\infty$
implies $(1 + \varepsilon f)^{-1} = ab \in H^\infty$. Hence
$f_{\varepsilon} = (\varepsilon I + f)(I+\varepsilon f)^{-1}
\in H^p$ for every $p \geq 1$, so if we can show that
$f_{\varepsilon} \in \M$, the proof is complete. That
$f_{\varepsilon} \in \M$ can be seen as follows:
$$
f_{\varepsilon} = \varepsilon (I + \varepsilon f)^{-1} + f(I+
\varepsilon f)^{-1};
$$
but $I = (I+\varepsilon f)(I+\varepsilon f)^{-1} = (I+ \varepsilon f)^{-1} +
\varepsilon f(I+\varepsilon f)^{-1}$ and $(I+\varepsilon f)^{-1}
\in \M$, so $\varepsilon f(I+\varepsilon f)^{-1} \in \M$,
implying $f(I+\varepsilon f)^{-1} \in \M$.

\smallskip
(3) $\Re\,(f_{\varepsilon}) = \Re\,\left((\varepsilon I + f)(I +
\varepsilon f)^{-1} \right) = \Re\, ([(\varepsilon I + \varepsilon^2 f) + (1
- \varepsilon^2) f)](I + \varepsilon f)^{-1}) = \varepsilon
I + (1 - \varepsilon^2) \Re\,(f(I + \varepsilon f)^{-1})$. Since we
assume that $\varepsilon < 1$, it is enough to show \noindent that $\Re\,(f(I+
\varepsilon f)^{-1}) \geq 0$. For this
\begin{align*}
\Re\, (f(I+\varepsilon f)^{-1}) &= \tfrac{1}{2}\left(f(I+\varepsilon
f)^{-1} + (I + \varepsilon f^*)^{-1}f^*\right) \\&= \tfrac{1}{2}(I +
\varepsilon f^*)^{-1} \left((I+\varepsilon f^*)f +
f^*(I+\varepsilon f)\right)(I+\varepsilon f)^{-1}  \\
&=\tfrac{1}{2}(I+\varepsilon f^*)^{-1}\left(2 \Re\, (f) + 2 \varepsilon
|f|^2\right)(I+\varepsilon f)^{-1} \geq 0.
\end{align*}
\smallskip

(4) 
We have for every $\varepsilon > 0$,
\begin{align*}
f_{\varepsilon} - f &= (\varepsilon I + f)(I +
\varepsilon f)^{-1} - f \\&= \left((\varepsilon I + f) - f(I+
\varepsilon f)\right)(I+\varepsilon f)^{-1} \\&= \varepsilon(I+
f^2)(I+\varepsilon f)^{-1},
\end{align*}
so
$$
\mu_t(f_\varepsilon - f) \pprec \varepsilon \mu_t(I+
f^2) \mu_t ((I+\varepsilon f)^{-1}).
$$
Since $\| (I + \varepsilon f)^{-1}\| \leq 1$, we get
$\mu_t((I+ \varepsilon f)^{-1}) \leq 1$ for every $t>0$.
 Also $I + f^2 \in
L^p(\M,\T)$ for every $p > 1$, so $\| f_{\varepsilon} - f \|_p
\leq \varepsilon \| I + f^2 \|_p \to 0$ (as $ \varepsilon \to 0$).
The proof is complete. 
\end{proof}

\begin{proposition}
Let  $u \in \M$ with $ u \geq 0$, and set $f = u + i \tilde{u}$.
There exists a constant $K$ (independent of $u$) such that, for every
$s > 0$, $$ \T (\chi_{(s, \infty)} (|f|)) \leq K \frac{\| u
\|_1}{s}.$$
\end{proposition}
\begin{proof}
 We follow  (at least in spirit) the argument of Helson in \cite{HEL}
for the commutative case.

Let $u$ and $f$ be as in the statement of the proposition, and
fix  $0 < \varepsilon < 1$. Set $f_{\varepsilon} $ as in  Lemma~2.
For $ s \in (0, \infty)$ fixed, consider the following transformation
on $\{z; \Re\,(z) \geq 0 \}$:
$$
A_s(z) = 1 + \frac{z-s}{z+s} \quad \text{ for all } z \in \{ w, \Re\,(w)
\geq 0 \}.
$$
It can be checked that the part of the plane $\{ z; |z| \geq s
\}$ is mapped to the half disk $\{ w;\,\Re\,(w) \geq 1 \}$;
this fact is very crusial in the argument of \cite{HEL} for the
 commutative case. Although we are unable to verify this fact through
functional calculus, 
one  can still recapture its consequences  by taking  the trace in every step.

Note that $\sigma (f_{\varepsilon} )$ is a compact subset of $\{ z;
\Re\,(z) \geq \varepsilon \}$. By the analytic
 functional calculus for Banach algebras, %
$$
A_s (f_{\varepsilon}) = I + (f_{\varepsilon} - s
I)(f_{\varepsilon} +  sI)^{-1} \in H^\infty
$$
and therefore (since $A_s$ is analytic)
\begin{equation}
\eqlabel{*}
\Phi(A_s(f_{\varepsilon})) = A_s (\Phi(f_{\varepsilon})).
\end{equation}
\begin{claim}
$\Phi(f_{\varepsilon}) = \Phi(u)_\varepsilon $.
\end{claim}

In fact $f_\varepsilon = (\varepsilon I + f) (I+ \varepsilon f)^{-1} \in
H^\infty$ and $ f_\varepsilon \cdot (I + \varepsilon f) = \varepsilon I +
f$,
so $ \Phi(f_{\varepsilon}) \Phi(I + \varepsilon f) = \Phi
(\varepsilon I + f)$. But $ \Phi(f) = \Phi(u)$, so we get $
\Phi(f_{\varepsilon}) (I + \varepsilon \Phi(u)) = \varepsilon I +
\Phi(u)$, and the claim follows.

\begin{claim}
$\Re \, (I+ (f_{\varepsilon} - sI)(f_{\varepsilon} + sI)^{-1}) \geq
0$.
\end{claim}

For this we have
\begin{align*}
\Re\, (I &+ (f_{\varepsilon} - sI)(f_{\varepsilon} + sI)^{-1}) = I
+ \tfrac{1}{2} \left((f_{\varepsilon} - sI)(f_{\varepsilon} + sI)^{-1}
+ (f^*_\varepsilon + sI)^{-1}(f^*_\varepsilon - sI)\right) \\
&= I + \tfrac{1}{2}(f^*_\varepsilon + sI)^{-1} \left((f^*_\varepsilon
+ sI)(f_{\varepsilon} - sI) + (f^*_\varepsilon -
sI)(f_{\varepsilon} + sI)\right)(f_{\varepsilon} + sI)^{-1} \\
&= I + (f^*_\varepsilon + sI)^{-1} \left(|f_{\varepsilon}|^2 -
s^2I\right)(f_{\varepsilon} + sI)^{-1} \\
&= (f^*_\varepsilon + sI)^{-1} \left((f^*_\varepsilon +
sI)(f_{\varepsilon} + sI)+ |f_{\varepsilon}|^2 - s^2I\right)
(f_{\varepsilon} + sI)^{-1} \\
& = 2 (f^*_\varepsilon + sI)^{-1} \left(|f_{\varepsilon}|^2 + s
\Re\,f_{\varepsilon}\right)(f_{\varepsilon} + sI)^{-1},
\end{align*}
and the claim follows from the fact that $\Re\,(f_{\varepsilon}) \geq
\varepsilon I$. 

Note that since $\Phi(u)$ is self-adjoint, so are
$\Phi(u)_\varepsilon$ and $A_s(\Phi(u)_\varepsilon)$.
 We conclude  from \eqref{*}
that $\T\left(I+ (\Phi(u)_\varepsilon - sI)
(\Phi(u)_\varepsilon + sI)^{-1}\right) \in \R$, and therefore
\begin{equation}\eqlabel{1}
\T\left(I+
\Re\,((f_{\varepsilon} - sI)(f_{\varepsilon} + sI)^{-1})\right) =
\T(A_s(\Phi(u)_\varepsilon)). 
\end{equation}

To prove the proposition, let $P=\chi_{(s,\infty)}(|f_{\varepsilon}|)$.
 The projection $P$ commutes with $|f_{\varepsilon}|$ and we  have
$$
\Re \left[I + (f_{\varepsilon} - sI)(f_{\varepsilon} + sI)^{-1}\right] =
(f^*_\varepsilon + sI)^{-1} \left[2|f_{\varepsilon}|^2 + 2 s
\Re\, (f_{\varepsilon})\right] (f_{\varepsilon} + sI)^{-1};
$$
but since $\Re\,(f_{\varepsilon}) \geq \varepsilon I$, we get
$$
2|f_{\varepsilon}|^2 + 2s \Re\,(f_{\varepsilon}) \geq
2|f_{\varepsilon}|^2 + 2s \varepsilon I,
$$
and hence %
\begin{equation}
\eqlabel{**}
\T\left[\Re\,(I + (f_{\varepsilon} - sI)(f_{\varepsilon}+sI)^{-1})\right] \geq
\T \left[(2|f_{\varepsilon}|^2 + 2s \varepsilon I)(f_{\varepsilon} + sI)^{-1}
(f_{\varepsilon}^* + sI)^{-1}\right]. 
\end{equation}
\begin{lemma}
Let  $a$ and $b$ be operators in 
 $\overline{ \M}$ with $ a \geq 0$, $b \geq 0$, and let $P$ be a
projection that commutes with $a$. Then $\T(ab) \geq \T(P(ab)P)$.
\end{lemma}
To see this, notice that, since $P$ commutes with $a$, $PaP \leq a$,
so $b^{1/2}PaPb^{1/2} \leq b^{1/2}ab^{1/2}$, implying that
$\T(b^{1/2}PaPb^{1/2} ) \leq \T(b^{1/2}ab^{1/2})$ and 
\begin{align*}
\T(P(ab)P) &= \T(P(ab)) \\
&= \T(PaPb) = \T(b^{1/2}PaPb^{1/2}) \\
&\leq \T(b^{1/2}a b^{1/2}) = \T(ab).
\end{align*}
The lemma is proved.

Applying Lemma~3 for $a = 2|f_{\varepsilon}| + 2s \varepsilon I$ and
$b = (f_{\varepsilon} + sI)^{-1} (f^*_\varepsilon + s
I)^{-1}$, we obtain %
\begin{equation}
\eqlabel{***}
\T\left[\Re\,(I + (f_{\varepsilon} -sI) (f_{\varepsilon} + sI)^{-1})\right]
\geq \T [ (2P |f_{\varepsilon}|^2 + 2 s \varepsilon P
)(f_{\varepsilon} + sI)^{-1} (f^*_\varepsilon + sI)^{-1}].
\end{equation}

Note that $(f_{\varepsilon} + sI)^{-1}(f^*_\varepsilon +
sI)^{-1} = \left(|f_{\varepsilon}|^2 + 2s \Re\,(f_{\varepsilon}) +
s^2I \right)^{-1}$ and $|f_{\varepsilon}|^2 + 2s \Re\, (f_{\varepsilon})
+ s^2I \leq |f_{\varepsilon}|^2 + 2s|f_{\varepsilon}| + s^2I$.

\begin{lemma}
Let $A,B,C$ be positive operators such that 
\begin{itemize}
\item[(i)] $A^{-1}$ and $B^{-1}$ exists, and
\item[(ii)] $A \leq B$.
\end{itemize}
Then $\T(C B^{-1}) \leq \T(C A^{-1})$.
\end{lemma}

To prove this lemma, observe that $\T (C B^{-1}) = \T(\alpha \beta )$,
where $\alpha = C^{1/2} B^{-1} A^{1/2}$ and $ \beta = A^{-1/2}
C^{-1/2}$. By H\"older's inequality, %
\begin{align*}
\T (C B^{-1}) &\leq \T( |\alpha|^2)^{1/2} \, \T(|\beta|^2)^{1/2} =
\T(\alpha^* \alpha)^{1/2}\, \T(\beta^* \beta)^{1/2} = \T(\alpha
\alpha^*)^{1/2}  \T(\beta^* \beta)^{1/2} \\
&= \T(C^{1/2} B^{-1} A^{1/2} A^{1/2} B^{-1} C^{1/2} )^{1/2}  
\, \T(C^{1/2} A^{-1/2} A^{-1/2} C^{1/2} )^{1/2} \\
&= \T(C^{1/2} B^{-1} A
B^{-1} C^{1/2})^{1/2}\,  \T(C A^{-1} )^{1/2}.
\end{align*}
But since $A \leq B$, we get $C^{1/2} B^{-1} (A) B^{-1} C^{1/2} \leq
C^{1/2} B^{-1} C^{1/2}$, 
and therefore
 $$\T(C B^{-1}) \leq \T(C B^{-1})^{1/2}\, \T(C
A^{-1})^{1/2}$$
 which shows that  $\T(C B^{-1}) \leq \T (C A^{-1}).$
The proof of Lemma~4 is complete.

Applying  Lemma~4 to
\begin{align*}
A &= (f^*_\varepsilon + sI) (f_{\varepsilon} + sI), \\
B &= |f_{\varepsilon}|^2 + 2s |f_{\varepsilon}|+ s^2I \\
\intertext{ and }
C &= 2|f_{\varepsilon}|^2 P + 2 s \varepsilon P,
\end{align*}
We obtain from \eqref{***} that
$$
\T \left[\Re\, (I + (f_{\varepsilon} - sI) (f_{\varepsilon} + sI)^{-1})\right]
\geq \T (C A^{-1}) \geq \T (C B^{-1})
$$
and hence %
$$
\T \left[ \Re\, (I + (f_{\varepsilon} - sI)(f_{\varepsilon} + sI)^{-1})\right]
\geq \T \left[(2|f_{\varepsilon}|^2P + 2s \varepsilon
P)(|f_{\varepsilon}|^2 + s|f_{\varepsilon}| + s^2I)^{-1}\right].
$$
Since $|f_{\varepsilon}|^2P \geq s^2P$, we get that 
\begin{equation}
\eqlabel{****}
\T\left[\Re\,(I + (f_{\varepsilon} - sI)(f_{\varepsilon} + sI)^{-1})\right]
 \geq \T
\left[P(|f_{\varepsilon}|^2 + 2 s \varepsilon I + s^2
I)(|f_{\varepsilon}|^2 + 2s |f_{\varepsilon}| + s^2 I)^{-1}\right].
\end{equation}
If we denote by $E^{|f_{\varepsilon}|}$ the spectral decomposition
of $|f_{\varepsilon}|$, then
$$
P(|f_{\varepsilon}|^2  + 2s \varepsilon I + s^2
I)(|f_{\varepsilon}|^2 + 2 s |f_{\varepsilon}| + s^2 I )^{-1}
= \int^{\infty}_s \frac{t^2 + 2 s \varepsilon + s^2}{t^2 + 2st + s^2}\, d
E^{|f_{\varepsilon}|}_t.
$$
Let 
$$\psi_{\varepsilon,s}(t) = \frac{t^2 + 2s \varepsilon + s^2}{t^2 + 2 st
+ s^2} \quad \text{ for } t \in [s, \infty).
$$ 
One can show that $\psi_{\varepsilon, s} $ attains its (unique) %
minimum at $t_0 = \varepsilon + \sqrt{\varepsilon^2 + \varepsilon s + s^2}$,
and therefore that
$$
\int^\infty_s \frac{t^2 + 2 s \varepsilon + s^2}{t^2 + 2 s t + s^2}\, d
E^{|f_{\varepsilon}|}_t \geq \psi_{\varepsilon, s}(t_0) P,
$$
so we deduce from \eqref{****} that
$$
\T\left[\Re\,(I + (f_{\varepsilon} - sI)(f_{\varepsilon} +
sI)^{-1})\right]\geq \psi_{\varepsilon, s}(t_0)\T(P).
$$

To finish the proof, recall from \eqref{1} that
$$
\T\left[\Re\,(I+(f_{\varepsilon}-sI)(f_{\varepsilon}+sI)^{-1})\right]=
\T\left[I+(\Phi(u)_\varepsilon - sI)(\Phi(u)_\varepsilon + sI)^{-1}\right],
$$
so
\begin{align*}
\T(P) &\leq \frac{1}{\psi_{\varepsilon, s}(t_0)}
 \T\left[I + (\Phi(u)_\varepsilon
- sI) (\Phi(u)_\varepsilon + sI)^{-1}\right] \\
&= \frac{1}{\psi_{\varepsilon, s}(t_0)} \T\left[2 \Phi(u)_\varepsilon
(\Phi(u)_\varepsilon + sI)^{-1}\right].
\end{align*}
But $(\Phi(u)_\varepsilon + sI)^{-1} =
 \frac{1}{s} (\frac{\Phi(u)_\varepsilon}{s} +
I)^{-1} $ has norm $\leq 1/s$, hence%
$$
\T(P) \leq \frac{2}{\psi_{\varepsilon, s} (t_0)} \frac{\| u_\varepsilon \|_1
}{s}.
$$
Now taking  $\varepsilon \to 0$, we get $\| u_\varepsilon \|_1 \to \| u \|_1$,
and $\T(P) \to \T (\chi_{(s, \infty)} (|f|))$. Note that 
$$
\psi_{\varepsilon, s} (t_0) = \frac{2 \varepsilon^2 + 3 \varepsilon s + 2 s^2 +
2 \varepsilon\sqrt{\varepsilon^2 + \varepsilon s + s^2}}{2 \varepsilon^2 + 3
\varepsilon s + 2 s^2 + 2 s (1 + \varepsilon)\sqrt{\varepsilon^2 +
 \varepsilon s + s^2}}
$$
so $ \lim_{\varepsilon \to 0} \psi_{\varepsilon, s} (t_0) = 1/2$.

Hence $\T(\chi_{(s, \infty)} (|f|)) \leq 4\, {\| u \|_1 }/s$.
The proof is complete.
\end{proof}

\

Recall that $L^{1, \infty}(\M,\T) = \{ a \in \overline{\M};\  \sup_{t > 0} t
\mu_t (a) < \infty \}$.

 Set $\|a\|_{1, \infty } = \sup_{t>0}t
\mu_t(a)$ for $a \in L^{1, \infty}(\M, \T)$. As in the commutative case,
$\| . \|_{1, \infty}$ is equivalent to a quasinorm in $L^{1, \infty}
(\M, \T)$, so there is a fixed constant $C$ such that, for every $a, b
\in L^{1, \infty} (\M, \T)$, we have $\|a + b \|_{1, \infty} \leq C(\| a
\|_{1, \infty} + \| b \|_{1, \infty})$. 

For $u \in \M$, let $Tu = u + i \tilde{u}$. From Theorem 1, $T$ is
linear and Proposition 1 can be restated as follows: 
$$
\text{ For any } u \in \M \text{ with } u \geq 0, \text{ we have } \|
Tu \|_{1, \infty} \leq 4 \| u \|_1;
$$
this implies that for $u \geq 0$,
$$
\| \tilde{u} \|_{1, \infty} \leq C ( 4+1) \| u \|_1 = 5C \| u \|_1.
$$
Now suppose  that $u \in \M$, $u = u^*$, $u = u_+ - u_-$ and $\tilde{u}
= \tilde{u}_+ - \tilde{u}_-$. Then
$$
\| \tilde{u} \|_{1, \infty} \leq C( \| \tilde{u}_+ \|_{1,\infty} + 
\| \tilde{u}_- \|_{1,\infty} ) \leq 5C^2 \| u \|_1.
$$
Similarly, if we require only $u \in \M$, %
we have $u = \Re\,(u)+ i \Im\, (u)$ and by linearity,
 $\tilde{u} = \widetilde{\Re(u) } + i\,
\widetilde{\Im(u)}$, and as above,
$$
\| \tilde{u} \|_{1, \infty} \leq 10C^3 \| u \|_1.
$$
We are now ready to extend ``$\sim$'' in $L^1(\M,\T)$: If $u \in
L^1(\M,\T)$, let $(u_n)_{n \in  \N}$ be a sequence in $\M$ such that $\| u -
u_n \|_1 \to 0$ as $n \to \infty$. Then
$$
\| \tilde{u}_n - \tilde{u}_m \|_{1, \infty} \leq 10 C^3 \| u_n - u_m \|_1,
$$
and since $\| u_n - u_m \|_{1} \to 0$ as $n,m \to \infty$, the
sequence $(\tilde{u}_n)_n$ converges in $L^{1, \infty}(\M, \T)$ to an
operator $\tilde{u}$. This defines $\tilde{u}$ for $u \in L^1(\M,\T) $.
 This
definition can be easily checked to be independent of the sequence
$(\tilde{u}_n)_n$ and agree with the conjugation operator  defined for
$p >1 $. 

 Letting $n \to \infty$ in the inequality $\| \tilde{u}_n
\|_{1, \infty} \leq 10 C^3 \| u_n \|_1$, we obtain the following
theorem ($H^{1,\infty}$ denotes the closure of $H^\infty$ in
 $L^{1,\infty}(\M,\T)$):
\begin{theorem}
There is a unique extension of ``$\sim$'' from $L^1(\M,\T) $
into $L^{1, \infty}(\M,\T)$ with the following property: $u + i \tilde{u}
\in H^{1, \infty}$ for all $ u \in L^1(\M,\T)$, and there is a constant
$K$ such that $\| \tilde{u} \|_{1, \infty } \leq K \| u \|_1$ for all
$ u \in L^1(\M,\T)$.
\end{theorem}

\begin{corollary}
For any $p$ with $0 < p < 1$ there exists a constant $K_p$ such that
$$\|\tilde{u}\|_p \leq K_p \|u\|_1 \ \ \text{for all}\  u \in L^1(\M,\T).$$
\end{corollary}

\begin{proof}
It is enough to show that such a constant exists for $u \in \M$, $u \geq 0$.
Recall that for $u \in \M$, the distribution $\lambda_s(u)$ equals
$\T(\chi_{(s, \infty)} (u))$.

 Let $F(s) = 1 - \lambda_s(u) = \T
(\chi_{(0,s)}(u))$. Assume that $\|u\|_1 \leq 1$. From proposition 1, 
$$
1- F(s) \leq \frac{4}{s}\|u\|_1 \leq \frac{4}{s}.
$$
Note that $F$ is a non-increasing right continuous function and for $p
> 0$, %
$$
\T (|\tilde{u}|^p) = \int^1_0 \mu_t(|\tilde{u}|)^p d t =
\int^\infty_0 s^p d F(s) \leq 1 + \int^\infty_1 s^p d F(s).
$$
If $A$ is
a point of continuity for $F ( A > 1)$, then 
$$
\int^A_1 s^p d F(s) = [s^p(F(s) - 1)]^A_1 + p\int^A_1 (1- F(s))s^{p-1}ds.
$$
Since $1 - F(s) \leq \frac{4}{s}$, we get that both $[s^p(F(s)
-1)]^A_1$ and $\int^A_1 (1 - F(s))s^{p-1} d s$ are bounded for $0 < p
< 1$, that is, $\int^1_0 \mu_t(|\tilde{u}|)^p dt$ has bound
independent of $u$.
\end{proof}

The {\bf Riesz projection} $\cal{R}$ can now be defined as
 in the commutative case:
for every $ a \in L^p(\M,\T)$, \ $(1\leq p \leq \infty)$,
   $$\cal{R}(a) = \frac{1}{2}(a + i\tilde{a} + \Phi(a)).$$

From Theorem~2, one can easily verify that $\cal{R}$ is
 a bounded projection from $L^p(\M,\T)$ onto $H^p$ for $1<p<\infty$.
 In particular $H^p$ is a complemented subspace of $L^p(\M,\T)$. For
$p=1$, Theorem~3 shows that $\cal{R}$ is bounded from $L^1(\M,\T)$  into
$H^{1,\infty}$.

Our next result gives a sufficient condition on an operator $a\in L^1(\M,\T)$
so that its conjugate $\tilde{a} \in L^1(\M,\T)$.
\begin{theorem}
 There exists a constant $K$ such that  for every $a \in \M$,
 $$ \| \tilde{a}\|_1 \leq K\T( |a|\log^+|a|) + K.$$
\end{theorem}

\begin{proof}
Since our proof of Theorem~2 follows the same line of argument as 
in \cite{DEV}, one can deduce as in \cite{HR}(Corollary~{2h}) that there 
is an absolute  constant $C$ such that
 $\| \tilde{a}\|_p \leq Cpq \|a\|_p$ for all
$a \in L^p(\M,\T)$, $1<p<\infty$ and $1/p +1/q=1$.
 The conclusion of the theorem
can be deduced as a straightforward adjustment of the commutative case 
in \cite{ZYG}(vol~{.II}, p.~{119}); we will present it here for completeness.
 
  Let $a \in \M$; we will assume first
 that $a\geq 0$. Let $(e_t)_t $ be the spectral 
decomposition of $a$. For each $k \in \N$, let $P_k =\chi_{[2^{k-1},2^k)}(a)$
be the  spectral projection relative to $[2^{k-1}, 2^k)$. Define $a_k =aP_k$
 for
$k \geq 1$ and $a_0=a\chi_{[0, 1)}(a)$.
 Clearly
$a= \sum_{k=0}^{\infty} a_k $ in $L^p(\M,\T)$ for every $1<p<\infty$.
 By linearity, $\tilde{a}=\sum_{k=0}^\infty \tilde{a}_k$. For every $k\in\N$,
 $\|\tilde{a}_k\|_1 \leq \|\tilde{a}_k\|_p \leq Cp^2{(p-1)}^{-1} \|a_k\|_p$.
 Since $a_k \leq 2^k P_k$, we get for $1<p<2$,
 $$\|\tilde{a}_k\|_1 \leq 4C \frac{1}{p-1}2^k \T(P_k)^{\frac{1}{p}}.$$
If we set $p =1 +\frac{1}{k+1}$ and $\epsilon_k =\T(P_k)$, we have
 $$ \|\tilde{a}_k\|_1 \leq 4C (k+1)2^k \epsilon_{k}^{\frac{k+1}{k+2}}.$$
Taking the summation over $k$,
 $$\|\tilde{a}\|_1 \leq \sum_{k=0}^\infty 4C(k+1)2^k \epsilon_{k}^{\frac{k+1}{k+2}}.$$
We note as in \cite{ZYG} that if $J=\{ k \in \N;\ \epsilon_k \leq 3^{-k}\}$ 
then 
 $$\sum_{k \in J} 4C(k+1)2^k \epsilon_{k}^{\frac{k+1}{k+2}} \leq
 \sum_{k=0}^\infty 4C(k+1)2^k (3^{-k})^{\frac{k+1}{k+2}} = \alpha <\infty.$$
On  the other hand, for $k \in \N \setminus J$, 
$\epsilon_{k}^{\frac{k+1}{k+2}} \leq \epsilon_k. 3^{\frac{k}{k+2}} \leq 
\beta \epsilon_k$ where $\beta=\sup_k 3^{\frac{k}{k+2}}$. So we get 
\begin{align*}
\|\tilde{a}\|_1 &\leq
 \alpha + 4C\beta \sum_{k=0}^\infty (k+1)2^k \epsilon_k \\
 &\leq \alpha +4C\beta(\epsilon_0 + 4\epsilon_1) +4C\beta \sum_{k\geq 2}(k+1)
2^k \epsilon_k.
\end{align*}
Since for $k\geq 2$,\   $k+1 \leq 3(k+1)$, we get
$$\|\tilde{a}\|_1 \leq 
\alpha +16C\beta +24C\beta\sum_{k\geq2}(k-1)2^{k-1}\epsilon_k.$$
To complete the proof, notice that  for $k\geq 2$,
\begin{align*}
 (k-1)2^{k-1} \epsilon_k &= \int_{2^{k-1}}^{2^k}  (k-1)2^{k-1}\ d\T(e_t)\\
  &\leq \int_{2^{k-1}}^{2^k} \frac{t\log t}{\log 2} \ d\T(e_t).
\end{align*}
Hence by setting $K_1=\max\{\alpha +16C\beta, 24C\beta/\log 2\}$, we get:
 $$\|\tilde{a}\|_1 \leq K_1 + K_1 \tau\left(a \log^+(a)\right).$$
 Now for a general element  $b \in
\M$, we decompose $b$ as $b=(b_{+}^{(1)} -b_{-}^{(1)}) + i( b_{+}^{(2)}-b_{-}^{(2)})$  where $b_{+}^{(i)}$ and $b_{-}^{(i)}$ are positive operators for
 $i=1, 2$.
One can easily verify that
 $\|\tilde{b}\|_1 \leq K + K\T(|b|\log^+(|b|))$ where $K = 4K_1$.
 The proof is complete.
\end{proof}
\begin{remark}
From Theorem~4, one can deduce that if $a \in L^1(\M,\T)$ is such that
$|a|\log^+(|a|)$ belongs to $L^1(\M,\T)$ then $\tilde{a} \in L^1(\M,\T)$.
\end{remark}

\section{Remarks on Szeg\"o's theorem}

Szeg\"o's theorem plays a very important role in theory of weak*-Dirichlet
algebras. It is still unknown if it has an extension in the context
of subdiagonal algebras.  In this section, we discuss different forms of
 possible extensions of Szeg\"o's theorem.  

\begin{proposition}
If $a$ and $a^{-1}$ belong to $H^\infty$, then
$$\inf\left\{\T(|a^*|^2 |I-f|^2);\ f \in H_{0}^\infty \right\}
 =\T(|\Phi(a)|^2).$$
\end{proposition}
\begin{proof}
Let $b=\Phi(a)$ and $p=I-ba^{-1}$. We will equip $\M$ with the following
 scalar product by setting for every $x, y \in \M$,
 $$\langle x,y\rangle =\T(|a^*|^2 y^*x).$$
One can easily verify, since $a$ is invertible, that $\M$ with
 $\langle.\,,.\rangle$ is a pre-Hilbertian. We denote the completion of
this space by $L^2(\M,|a^*|^2)$. Let $B$ be the closure of $H_{0}^\infty$ in
the space $L^2(\M,|a^*|^2)$. We claim that $p$ is the projection of $I$
 into $B$.

  First we will show that $p \in H_{0}^\infty$ (and thus  $p \in B$):
clearly $p \in H^\infty$; also
$\Phi(p)= I-b\Phi(a^{-1})$ and since both $a$ and $a^{-1}$ belong to
 $H^\infty$, $I=\Phi(aa^{-1})=\Phi(a).\Phi(a^{-1})$ so $\Phi(a^{-1})=b^{-1}$
which implies that $\Phi(p)=0$.

 To prove the claim it is enough
to check that $I-p=ba^{-1} \perp H_{0}^\infty$ for $\langle.\,,.\rangle$:
for $f \in H_{0}^\infty$, we have
\begin{align*}
 \langle f, I-p \rangle &=\T(|a^*|^2 (ba^{-1})^*f) \\  
                        &=\T(aa^*((a^{-1})^*b^*)f)\\
                        &=\T(ab^*f)
\end{align*}
and since $b \in D$, $ab^* \in H^\infty$ so $\T(ab^*f)=0$. 

To complete the proof of the theorem, note that  
$\text{dist}(I, B)= || I-p ||_{L^2(|a^*|^2)}$ so 
 $\inf\left\{||I-f||_{L^2(|a^*|^2)};\ f \in B \right\}=||I-P||_{L^2(|a^*|^2)}$
i.e 
$$\inf\left\{\T(|a^*|^2 |I-f|^2);\ f \in H_{0}^\infty \right\}=
\T(|a^*|^2 |ba^{-1}|^2)=\T(|b|^2).$$
The proof is complete.
\end{proof}
\begin{remark}  Using similar  argument as above with $p=I-a^{-1}b$ and
the Hilbert space defined by
$\langle x,y\rangle = \T(|a|^2 xy^*)$, one can deduce the following identity:
 $$\inf\left\{\T(|a|^2 |I-f^*|^2);\ f \in H_{0}^\infty\right\}=
\T(|\Phi(a)|^2).$$
\end{remark}

We remark also that if $v \in \M$ is such that  $v=v^*$ then there exists
 $a \in H^\infty$ with $a^{-1} \in H^\infty$ satisfying $|a^*|^2= e^v$; in
fact for such $v$, $e^{v/2}$ is invertible in $\M$. Apply \cite{SAI1}
(Proposition~1) to get operators $a \in H^\infty$ and $u$ unitary in $\M$ 
such that $au=e^{v/2}$. Since $u$ and $e^{v/2}$ are invertible, it is clear
that $a$ is invertible. The fact that $|a^*|^2 = e^v$ is immediate. Similarly,
an operator $b \in H^\infty$ with $b^{-1} \in H^\infty$ can also be choosen so that $|b|^2=e^{v}$.
  
\begin{proposition}
Let $h \in \M$, $h \geq 0$ and $\log(h)$ exists then
 $$\exp(\T(\log(h))) \geq \inf\left\{ \T(h e^{\Re{f}});\ f \in \M,
\T(f)=0\right\}.$$
\end{proposition}
The proof of Proposition~3 is based on the following simple extension
 of the usual Jensen's inequality:
\begin{lemma}
Let $h \in\M$, $h\geq 0$ then 
$\T(\log(h)) \leq \log(\T(h)).$
\end{lemma}
\begin{proof}
Let $\varphi: [0, \infty) \to \R$ defined by $\varphi(x)=\log(x+1)$. The function $\varphi$ is continuous, increasing and $\varphi(0)=0$. We get
from \cite{FK}(Corollary~2.8) and the usual Jensen's inequality for functions
 that  
\begin{align*}
 \T(\log(h +I)) &=\int_{0}^{1} \mu_t( \log(h) + I) \ dt \\
   &=\int_{0}^{1} \log(\mu_t(h) +1)\ dt \\
   &\leq \log\left( \int_{0}^{1} \mu_t(h)\ dt + 1 \right) \\
   &=\log\left( \T(h + I)\right).
\end{align*}
So we have $\T(\log(h +I)) \leq \log( \T(h +I))$.  Fix $\varepsilon >0$.
Applying the same inequality for $h/\varepsilon$,  we get
\begin{align*}
 \T(\log(h/\varepsilon +I) &\leq \log\left( \T(h/\varepsilon) +1\right) \\
  \T(\log(h +\varepsilon I) -(\log\varepsilon)I) &\leq
 \log\left( \T(h + \varepsilon I)\right)  - \log\varepsilon \\
 \T(\log(h +\varepsilon I)) &\leq \log\left(\T(h)  +\varepsilon \right).
\end{align*}
By letting $\varepsilon \to 0$, the desired inequality follows.
 The lemma is proved.
\end{proof}

\noindent
{\it Proof of Proposition~3.}
Let $g \in \M$ such that  $g=g^*$, $g$ commutes with $h$ and $\T(g)=0$.
Applying  Lemma~5 to the operator $he^{g}$, we have
$\exp\left( \log(h) \right) \leq \T(he^{g})$ and therefore,
 $$\exp\left( \T(\log(h))\right) \leq
 \inf\left\{ \T(he^{g});\ g=g^*, g\  \text{commutes with}\  h,
 \T(g)=0\right\}.$$
Let $\lambda= \T(\log(h))$ and $g=\lambda I - \log(h)$. Clearly
$g=g^*$, $g$ commutes with $h$ and $\T(g)=0$ and it is easy to check that
$\T(he^{g})= \exp\left(\T(\log(h))\right)$
 so the inequality above is in fact an equality i.e
 $$\exp\left( \T(\log(h))\right)
=\inf\left\{ \T(h e^g);\ g=g^*, g\ \text{commutes with}\ h, \T(g)=0 \right\}.$$
This implies that
$$\exp\left( \T(\log(h))\right)
 \geq \inf\left\{ \T(h e^g);\ g=g^*,  \T(g)=0 \right\}.$$
The proof is complete.
\qed

The above proposition leads to the following question:
If $h \in H^\infty$, is it true that $\exp(\T(\log|h|))\geq |\T(h)|$ ?
This inequality is known as Jensen's inequality for $H^\infty$. 
\begin{remark}
In \cite{DA}, characterization of real functions in $L^1(\mathbb{T})$ that 
have rearrangemment in $\Re{H_{0}^1(\mathbb{T})}$ 
 were given (see also \cite{KA2} for another proof).
 The same characterization was shown to be true 
for the weak*-Dirichlet algebra setting in \cite{Lan2}. It would be
 interesting to know  if such characterization holds for the non-commutative
case. We note that the proofs given in \cite{KA2} and \cite{Lan2} 
made use  of Szeg\"o's theorem in a very crusial way.
\end{remark}

\end{document}